\documentclass[12pt]{amsart}
\usepackage{fullpage,url,amssymb}

\DeclareFontEncoding{OT2}{}{} 

\usepackage{color}

\newcommand{\defi}[1]{\textsf{#1}} 

\newenvironment{romanenum}{\hfill \begin{enumerate} }{\end{enumerate}}

\newcommand{\C}{{\mathbb C}}
\newcommand{\F}{{\mathbb F}}

\newcommand{\PP}{{\mathbb P}}
\newcommand{\Q}{{\mathbb Q}}

\newcommand{\Z}{{\mathbb Z}}
\newcommand{\Qbar}{{\overline{\Q}}}

\newcommand{\Fbar}{{\overline{\F}}}

\DeclareMathOperator{\Ig}{Ig}

\DeclareMathOperator{\Char}{char}

\DeclareMathOperator{\Aut}{Aut}
\DeclareMathOperator{\Gal}{Gal}

\newcommand{\tors}{{\operatorname{tors}}}

\newcommand{\GL}{\operatorname{GL}}
\newcommand{\SL}{\operatorname{SL}}

\newcommand{\PSL}{\operatorname{PSL}}

\newcommand{\surjects}{\twoheadrightarrow}

\newcommand{\isom}{\simeq}

\newcommand{\intersect}{\cap} 

\newtheorem{theorem}{Theorem}[section]
\newtheorem{lemma}[theorem]{Lemma}
\newtheorem{corollary}[theorem]{Corollary}
\newtheorem{proposition}[theorem]{Proposition}

\theoremstyle{definition}

\theoremstyle{remark}
\newtheorem{remark}[theorem]{Remark}

\usepackage[
        backref,
        pdfauthor={Bjorn Poonen}, 
        pdftitle={Gonality of modular curves in characteristic p},
]{hyperref}
\usepackage[alphabetic,backrefs,lite]{amsrefs} 

\begin{document}

\title[Gonality of modular curves]{Gonality of modular curves in characteristic~$p$}
\subjclass[2000]{Primary 14G35; Secondary 14H51, 11G18}
\keywords{gonality, modular curve, Igusa curve, image of Galois}
\author{Bjorn Poonen}
\thanks{This research was supported by NSF grant DMS-0301280
          and the Miller Institute for Basic Research in Science.}
\address{Department of Mathematics, University of California, 
        Berkeley, CA 94720-3840, USA}
\email{poonen@math.berkeley.edu}
\urladdr{http://math.berkeley.edu/\~{}poonen}
\date{February 23, 2006}

\begin{abstract}
Let $k$ be an algebraically closed field of characteristic $p$.
Let $X(p^e;N)$ be the curve parameterizing
elliptic curves with full level $N$ structure 
(where $p \nmid N$) and full level $p^e$ Igusa structure.
By modular curve, we mean a quotient of any $X(p^e;N)$
by any subgroup of
$\left( (\Z/p^e\Z)^\times \times \SL_2(\Z/N\Z) \right)/\{\pm 1\}$.
We prove that in any sequence of distinct modular curves over $k$,
the $k$-gonality tends to infinity.
This extends earlier work, in which the result
was proved for particular sequences of modular curves,
such as $X_0(N)$ for $p \nmid N$.
We give an application to the function field analogue
of a uniform boundedness statement for the image of Galois
on torsion of elliptic curves.
\end{abstract}

\maketitle

\section{Introduction}\label{S:introduction}

\subsection{Gonality}

The \defi{gonality} $\gamma_k(X)$ of a curve\footnote{%
All our curves are assumed to be geometrically integral.%
}
 $X$ over a field $k$
is the smallest possible degree of a 
dominant rational map $X \dashrightarrow \PP^1_k$.
For any field extension $L$ of $k$, we define also the \defi{$L$-gonality}
$\gamma_L(X)$ of $X$ as the gonality 
of the base extension $X_L:=X \times_k L$.
It is an invariant of the function field $L(X)$ of $X_L$.

\subsection{General facts about gonality}
\label{S:facts}

The following facts are well-known.
\begin{proposition}
\label{P:facts}
Let $X$ be a curve of genus $g$ over a field $k$.
\begin{romanenum}
\item 
If $L$ is a field extension of $k$, then $\gamma_L(X) \le \gamma_k(X)$.
\item
If $k$ is algebraically closed, and $L$ is a field extension of $k$,
then $\gamma_L(X) = \gamma_k(X)$.
\item 
If $g>1$, then $\gamma_k(X) \le 2g-2$.
For each $g>1$, there exist $k$ and $X$ for which equality holds.
\item 
If $X(k) \ne \emptyset$, then $\gamma_k(X) \le g+1$.
If $X(k) \ne \emptyset$ and $g \ge 2$, then $\gamma_k(X) \le g$.
\item
If $k$ is algebraically closed,
then $\gamma_k(X) \le \lfloor \frac{g+3}{2} \rfloor$.  
Equality holds for a general curve of genus $g$ over $k$.
\item 
If $\pi\colon X \dashrightarrow Y$ is a dominant rational map 
of curves over $k$, 
then $\gamma_k(X) \le (\deg \pi) \gamma_k(Y)$.
\item
If $\pi\colon X \dashrightarrow Y$ is a dominant rational map 
of curves over $k$, then $\gamma_k(Y) \le \gamma_k(X)$.
\end{romanenum}
\end{proposition}

\begin{proof}\hfill
\begin{romanenum}
\item
Trivial.
\item
Given a map over $L$, 
standard specialization arguments give a map over $k$
of the same degree.
\item
The canonical linear system $|K|$ has dimension $g-1 \ge 1$
and degree $2g-2$.
So we may use a rational function whose divisor is the difference of 
two different canonical divisors.
Equality holds for the general curve 
over the function field of the moduli space of genus-$g$ curves
in characteristic $0$,
since its only line sheaves are the powers of the canonical sheaf:
this was the Franchetta conjecture, proved in \cite{Harer1983}
and strengthened in \cite{Mestrano1987}.
\item
Let $P \in X(k)$.
The Riemann-Roch theorem shows that $\dim |(g+1)P| \ge 1$,
and that $\dim |K -(g-2)P| \ge 1$ if $g-2 \ge 0$.
These linear systems have degree $g+1$ and $g$, respectively.
\item
These are consequences of Brill-Noether theory: see (1.1) and (1.5) 
of \cite{Arbarello-et-al1985}*{Chapter~V} for an exposition.
The first statement is proved in arbitrary characteristic 
in \cites{Kleiman-Laksov1972,Kleiman-Laksov1974}.
The second statement is proved in characteristic $0$
in \cites{Farkas1966,Martens1967,Martens1968},
and can be deduced in characteristic $p$ 
from the unramified case of \cite{Osserman2005}*{Theorem~1.2}, for instance.
\item
Trivial.
\item
(The ideas in the following argument 
go back at least to \cite{Newman1972}*{Theorem~VII.2}.)
Choose $f \in k(X)$ of degree $d:=\gamma_k(X)$.
Let $r=\deg \pi = [k(X):k(Y)]$.
Let $P(T) \in k(Y)[T]$ be the characteristic polynomial of $f$
viewed as an element of the field extension $k(X)$ of $k(Y)$.
For some finite normal extension $M$ of $k(Y)$,
we may write $P(T)=\prod_{i=1}^r (T-f_i)$ for some $f_i \in M$.
As a function in $M$, $f$ has degree $[M:k(X)] d$.
The same is true of each $f_i$, since they are all in the same
$\Aut(M/k(Y))$-orbit.
The polar divisor of a coefficient of $P$ viewed in $M$
is at most the sum of the polar divisors of the $f_i$,
so each coefficient
has degree at most $r [M:k(X)] d = [M:k(Y)] d$ as a function in $M$,
and hence degree at most $d$ as a function in $k(Y)$.
Since $f$ is non-constant, at least one of these coefficients
is non-constant.
Thus $\gamma_k(Y) \le d$.
\end{romanenum}
\end{proof}

\subsection{Modular curves}

Our goal is to obtain lower bounds on the gonality of modular curves.
By Proposition~\ref{P:facts}(ii), it suffices to consider
$k=\C$ and $k=\Fbar_p$ (an algebraic closure of $\F_p$) for each prime $p$.

Suppose $N$ is a positive integer not divisible by the characteristic of $k$.
Choose a primitive $N$-th root of unity $\zeta$.
Let $X(N)$ be the smooth projective model of the (possibly coarse) 
moduli space parameterizing triples $(E,P,Q)$ where $P,Q \in E$ 
are a basis for $E[N]$ with Weil pairing $e_N(P,Q)=\zeta$.
For $k=\C$, we can describe $X(N)(\C)$ alternatively as the quotient
of an extended upper half plane by a finite-index subgroup of $\PSL_2(\Z)$.
More generally, any congruence subgroup $G \le \PSL_2(\Z)$
gives rise to a curve $X_G$ over $\C$.
We have the following theorem of Abramovich:

\begin{theorem}[\cite{Abramovich1996}]
\label{T:Abramovich}
Let $D=(\PSL_2(\Z):G)$.
Then $\gamma_\C(X_G) \ge \frac{7}{800} D$.
\end{theorem}

\begin{remark}
As mentioned in \cite{Abramovich1996},
combining Theorem~\ref{T:Abramovich}
with the genus bound $g-1 \le D/12$ \cite{Shimura1994}*{Proposition~1.40}
yields
\[
        \gamma_\C(X_G) \ge \frac{21}{200}(g-1).
\]
The proof of Theorem~\ref{T:Abramovich} 
makes use of a lower bound on the leading nontrivial eigenvalue
of the noneuclidean Laplacian; this bound has been improved since 1996,
so the constants $7/800$ and $21/200$ can be improved too.
See also \cite{Baker-et-al2005}*{\S4.3} for some further results.
\end{remark}

In characteristic $p$, one has other kinds of modular curves,
involving level structure where $p$ divides the level.
If $q=p^e$ for some $e \in \Z_{\ge 1}$, the \defi{Igusa curve of level $q$} 
is the smooth projective model $\Ig(q)$ of the curve over $\F_p$
parameterizing pairs $(E,R)$ where $E$ is an ordinary elliptic curve,
and $R$ is a generator of the kernel of the degree-$q$ Verschiebung isogeny
$V_q \colon E^{(q)} \to E$, where $E^{(q)}$ is the elliptic curve obtained by
raising all the coefficients of a model of $E$ to the $q$-th power.
Given $N$ not divisible by $p$, and $q=p^e$,
we can define also a hybrid modular curve $X(p^e;N)$ over $\Fbar_p$
parameterizing $(E,R,P,Q)$, with $R \in \ker V_q$ and $P,Q \in E[N]$
as above.

The group $G_{p^eN}:=(\Z/p^e\Z)^\times \times \SL_2(\Z/N\Z)$ 
acts on $X(p^e;N)$.
The kernel of the action is $\{\pm 1\}$ embedded diagonally in $G_{p^e N}$.
For any subgroup $G \le G_{p^e N}$ containing $\{\pm 1\}$, let $X_G$
be the smooth projective model of the quotient $X(p^e;N)/G$.
The $G_{p^e N}$ form an inverse system with inverse limit
\[
        S:=\Z_p^\times \times \prod_{\text{prime $\ell \ne p$}} \SL_2(\Z_\ell).
\]
The inverse image of $G$ under $S \surjects G_{p^e N}$
is an open subgroup of the profinite group $S$,
and every open subgroup of $S$ containing $\{\pm 1\}$
arises this way for some $p^e N$.
Thus we may define $X_G$ for any open subgroup $G$ of $S$ 
containing $\{\pm 1\}$.

It seems likely that there is a constant $c>0$ independent 
of $p$ and $G$ such that $\gamma_{\Fbar_p}(X_G) \ge c (S:G)$.
We are unable to prove such a linear lower bound, even for fixed $p$,
but we can show that the gonality goes to infinity for fixed $p$.
Here is our main theorem:

\begin{theorem}
\label{T:main}
Fix a prime $p$.
Let $G_1,G_2,\ldots$ be a sequence of distinct open subgroups of $S$ 
containing $\{\pm 1\}$.
Then $\gamma_{\Fbar_p}(X_{G_i}) \to \infty$ as $i \to \infty$.
\end{theorem}

\subsection{Outline of proof of main theorem}

Many of the ideas used in the proof of Theorem~\ref{T:main} 
are due to earlier authors,
though we consider a broader class of modular curves
than had been treated earlier.
Section~\ref{S:change} proves 
Theorem~\ref{T:change},
an inequality in the direction opposite to Proposition~\ref{P:facts}(i):
the ideas used here and their application to 
the classical modular curves $X_0(N)$ can be found in
\cite{Harris-Silverman1991}, 
\cite{Nguyen-Saito1996preprint}, 
and \cite{Baker-thesis}*{Chapter 3}.
Theorem \ref{T:change}
reduces the problem to finding lower bounds on gonality over 
{\em finite} fields,
and these can be obtained by counting in the spirit of \cite{Ogg1974},
which among other things 
determined the $N$ for which $X_0(N)$ is hyperelliptic.
In Section~\ref{S:prime to p} we find that, as in \cite{Ogg1974}, 
modular curves of level prime to $p$ 
have too many supersingular points over $\F_{p^2}$
to have small gonality.
In Section~\ref{S:power of p} we cite results
of Schweizer \cite{Schweizer2005}, who obtained lower bounds
on the $\Fbar_p$-gonality of Igusa curves directly from the geometry
of the curves, instead of first getting lower bounds on $\F_p$-gonality
by counting $\F_p$-points.
Section~\ref{S:group theory} uses Goursat's lemma to study the subgroups
of $S$, so that the prime-to-$p$ and $p$-power cases can be combined
to prove the general case of Theorem~\ref{T:main} 
in Section~\ref{S:general case}.

\subsection{Application to the image of Galois}

One application of results like Theorem~\ref{T:main},
noted already by many other authors,
is to the function field analogue of
the strong uniform boundedness theorem for elliptic curves.
By the work of Mazur, Kamienny, and Merel \cite{Merel1996},
for every $d \in \Z_{\ge 1}$,
there exists a constant $N_d$ such that
for any number field $K$ with $[K:\Q] \le d$
and for any elliptic curve $E$ over $K$,
the torsion subgroup $E(K)_\tors$ of 
the finitely generated abelian group $E(K)$
satisfies $\#E(K)_\tors \le N_d$.
In the function field case, we can prove a stronger result,
one which bounds the index of the image of Galois acting on torsion.
If $E$ is an elliptic curve over a field $K$ of characteristic $p \ge 0$,
and $K^s$ is a separable closure of $K$,
there exists a homomorphism
\[
        \rho_E\colon \Gal(K^s/K) \to \Z_p^\times \times \prod_{\ell \ne p} \GL_2(\Z_\ell)
\]
describing the Galois action on the Tate modules of $E$.
(Of course, if $\Char K=0$, there is no $\Z_p^\times$ factor.)

\begin{theorem}
\label{T:uniform boundedness}
Given $p \ge 0$ and $d \in \Z_{\ge 1}$,
there exists a constant $N_{p,d}$ such that
for any field $k$ of characteristic $p$,
any field $K$ of degree $\le d$ over $k(t)$,
and any elliptic curve $E$ over $K$ with $j(E)$ not algebraic over $k$,
the index $(S:\rho_E(\Gal(K^s/K)) \intersect S)$ is at most $N_{p,d}$.
\end{theorem}

\begin{remark}
In Theorem~\ref{T:uniform boundedness},
we cannot hope to bound the index of $\rho_E(\Gal(K^s/K))$
in $\Z_p^\times \times \prod_{\ell \ne p} \GL_2(\Z_\ell)$
(i.e., with $\GL_2$ instead of the $\SL_2$ in the definition of $S$),
since the determinant of the image in $\GL_2(\Z_\ell)$
gives the action of $\Gal(K^s/K)$ on roots of unity,
and this is trivial if $k$ is algebraically closed, for example.
\end{remark}

Theorem~\ref{T:uniform boundedness} will be deduced from 
Theorem~\ref{T:main}
in Section~\ref{S:uniform boundedness}.

\section{Change in gonality under extension of the ground field}
\label{S:change}

In this section we give an exposition of the ``tower theorem'' of 
Nguyen and Saito \cite{Nguyen-Saito1996preprint}*{Theorem 2.1},
and its implication for relating gonalities of a single curve 
over different fields.
We will reprove it as our Proposition~\ref{P:tower}, 
since \cite{Nguyen-Saito1996preprint} remains unpublished after 10 years, 
and since we can simplify the proof slightly.
Throughout this section, $k$ is a perfect field.

\begin{proposition}[Castelnuovo-Severi inequality]
\label{P:Castelnuovo-Severi}
Let $F$, $F_1$, $F_2$ be function fields of curves over $k$,
of genera $g$, $g_1$, $g_2$, respectively.
Suppose that $F_i \subseteq F$ for $i=1,2$
and the compositum of $F_1$ and $F_2$ in $F$ equals $F$.
Let $d_i=[F:F_i]$ for $i=1,2$.
Then
\[
	g \le d_1 g_1 + d_2 g_2 + (d_1-1)(d_2-1).
\]
\end{proposition}

\begin{proof}
See \cite{Stichtenoth1993}*{III.10.3}.
\end{proof}

Let $X$ be a curve over $k$.
A subfield $F$ of $k(X)$ will be called \defi{$d$-controlled} if 
there exists $e \in \Z_{>0}$ such that
$[k(X):F]=d/e$ and the genus of $F$ is $\le (e-1)^2$.

\begin{lemma}
\label{L:d-controlled}
If $F$ is $d$-controlled, and $f \in k(X)$ is a rational function of
degree $d$, then $F(f)$ is $d$-controlled.
\end{lemma}

\begin{proof}
View $F(f)$ as the compositum of $F$ and $k(f)$ in $k(X)$.
Let $a=[F(f):F]$.
Then $[k(X):F(f)] = d/(ae)$, so $[F(f):k(f)]=ae$.
By Proposition~\ref{P:Castelnuovo-Severi}, the genus of $F(f)$ is at most
\[
	a (e-1)^2 + 0 + (a-1)(ae-1)
	\quad=\quad
	(ae-1)^2 - ae (a-1)(e-1) 
	\quad\le\quad 
	(ae-1)^2,
\]
since $a,e \ge 1$.
\end{proof}

\begin{corollary}
\label{C:d-controlled}
A subfield of $k(X)$ generated over $k$ by one or more elements of degree $d$
is $d$-controlled.
\end{corollary}

\begin{proof}
Induction on the number of elements:
the case of one element is trivial ($e=1$),
and Lemma~\ref{L:d-controlled} gives the inductive step.
\end{proof}

\begin{proposition}[Tower theorem]
\label{P:tower}
Let $X$ be a curve over a perfect field $k$.
Let $L \supseteq k$ be an algebraic field extension.
Let $d=\gamma_L(X)$.
Then $k(X)$ has a $d$-controlled subfield.
\end{proposition}

\begin{proof}
Enlarging $L$ cannot increase $\gamma_L(X)$, so we may assume $L/k$ is Galois.
Choose $f \in L(X)$ of degree $d$.
Let $F_L$ be the subfield generated over $L$ by the $\Gal(L/k)$-conjugates 
of $f$.
By Corollary~\ref{C:d-controlled}, $F_L$ is $d$-controlled
as a subfield of $L(X)$.
The action of $\Gal(L/k)$ on $L(X)$ preserves $F_L$,
and the invariant subfield $F_k:=F_L^{\Gal(L/k)}$
satisfies $[k(X):F_k]=[L(X):F_L]$ and has the same genus as $F_L$.
Thus $F_k$ is a $d$-controlled subfield of $k(X)$.
\end{proof}

\begin{theorem}
\label{T:change}
Let $X$ be a curve over a perfect field $k$.
Let $L \supseteq k$ be an algebraic field extension.
Let $d=\gamma_L(X)$.
Assume that $X(k) \ne \emptyset$.
\begin{romanenum}
\item
If $d \le 2$, then $\gamma_k(X)=d$.
\item \label{I:d>2}
If $d>2$, then $\gamma_k(X) \le (d-1)^2$.
\item \label{I:simple square root}
In any case, $\gamma_L(X) \ge \sqrt{\gamma_k(X)}$.
\end{romanenum}
\end{theorem}

\begin{proof}\hfill
\begin{romanenum}
\item
If $d=1$, then $X \isom \PP^1_k$, so $\gamma_k(X)=1$.
If $d=2$, then $X$ is elliptic or hyperelliptic;
if elliptic, then $\gamma_k(X)=2$;
if hyperelliptic then the canonical map is a degree-$2$ map to a
genus-$0$ curve $Z$ over $k$, and $Z(k) \ne \emptyset$ so $Z \isom \PP^1_k$,
so $\gamma_k(X)=2$.
\item
Now suppose $d>2$.
By Proposition~\ref{P:tower} there exists $e \in \Z_{>0}$ and a
rational map $\pi\colon X \dashrightarrow Y$ of curves over $k$
such that $\deg \pi = d/e$ and the genus $g$ of $Y$ satisfies $g \le (e-1)^2$.
We have $Y(k)\ne\emptyset$.
If $g=0$, then $Y \isom \PP^1$, so 
\[
	\gamma_k(X) \le d/e \le d < (d-1)^2.
\]
If $g=1$, then $e \ge 2$ and $\gamma_k(Y)=2$, so
\[
	\gamma_k(X) \le (d/e) \gamma_k(Y) \le (d/2) 2 = d < (d-1)^2.
\]
If $g \ge 2$, then $\gamma_k(Y) \le g$ by Proposition~\ref{P:facts}(iv), so
\[
	\gamma_k(X) \le \frac{d}{e} \gamma_k(Y) \le \frac{d}{e} (e-1)^2.
\]
For $e \in [1,d]$, the function $\frac{d}{e}(e-1)^2$ is maximized at $e=d$, 
and the value there is $(d-1)^2$.
\item
This follows directly from the first two parts.
\end{romanenum}
\end{proof}

\begin{remark}
The hypothesis $X(k)\ne \emptyset$ is necessary:
Genus-$1$ curves over $\Q$ have $\Qbar$-gonality $2$,
but their $\Q$-gonality can be arbitrarily large.
\end{remark}

\begin{remark}
We do not know whether the $(d-1)^2$ in Theorem~\ref{T:change}\eqref{I:d>2} 
can be improved.
\end{remark}

\begin{remark}
For $N=38, 44, 53, 61$, the modular curve $X_0(N)$ is of genus $4$
and has $\Q$-gonality $4$ and $\Qbar$-gonality $3$ 
\cite{Hasegawa-Shimura1999}*{p.~136}.
In particular, Theorem~\ref{T:change}\eqref{I:d>2} is best possible for $d=3$.
\end{remark}

\section{Level prime to $p$}
\label{S:prime to p}

Suppose $p \nmid N$.
We begin by defining a twisted form $X(N)'$ over $\F_{p^2}$ of $X(N)$.
Let $M$ be $(\Z/N \Z)^2$ made into a 
$\Gal(\Fbar_p/\F_{p^2})$-module
by letting the $p^2$-power Frobenius automorphism act as 
multiplication by $-p$.
There exists an isomorphism of Galois modules 
$\iota\colon \bigwedge^2 M \to \mu_N$; 
fix one.
Let $X(N)'$ be the smooth projective model of 
the affine curve over $\F_{p^2}$
parameterizing pairs $(E,\phi)$ where $E$ is an elliptic curve
and $\phi$ is an isomorphism $E[N] \to M$
under which the Weil pairing corresponds to $\iota$.
Over $\Fbar_p$, $X(N)'$ becomes isomorphic to $X(N)$.
The automorphisms of $M$ as an abelian group automatically commute with
the Galois action, so they induce automorphisms of $X(N)'$
defined over $\F_{p^2}$.
Thus we get 
$\SL_2(\Z/N\Z)/\{\pm 1\} \le \operatorname{Aut} X(N)'$.
Moreover, it follows from \cite{Baker-et-al2005}*{Lemma~3.21}
that all the points of $X(N)'$ corresponding to 
supersingular elliptic curves are defined over $\F_{p^2}$.

\begin{proposition}
\label{P:special modular}
Let $p$, $N$, and $X(N)'$ be as above.
Let $G$ be a subgroup of $\SL_2(\Z/N\Z)/\{\pm 1\}$ of index $D$.
Let $X$ be the curve $X(N)'/G$.
Then the $\F_{p^2}$-gonality $\gamma$ of $X$ satisfies
\[
        \gamma \ge \frac{p-1}{12(p^2+1)} D.
\]
\end{proposition}

\begin{proof}
This resembles the proof of \cite{Baker-et-al2005}*{Lemma~3.22}.
By \cite{Baker-et-al2005}*{Lemma~3.20},
the number of supersingular points on $X$ is $\ge (p-1)D/12$,
and these are images of supersingular points on $X(N)'$
so they are defined over $\F_{p^2}$;
thus $\#X(\F_{p^2}) \ge (p-1)D/12$.
On the other hand, 
$\#X(\F_{p^2}) \le \gamma \# \PP^1(\F_{p^2}) = \gamma (p^2+1)$.
Combine the two previous sentences.
\end{proof}

\begin{remark}
Let $g$ be the genus of $X$.
One could also combine Proposition~\ref{P:special modular} with
the bound $g-1 < D/12$ of \cite{Shimura1994}*{Proposition~1.40}
to give a lower bound for $\gamma$ in terms of $g$ instead of $D$.
\end{remark}

We now consider the $\Fbar_p$-gonality of 
all modular curves of level prime to $p$.

\begin{corollary}
\label{C:Phi}
Fix $p$.
Define $\Phi_p(D):=\sqrt{\frac{p-1}{12(p^2+1)} D}$.
If $G$ is the inverse image under 
$S \surjects \prod_{\ell \ne p} \SL_2(\Z_\ell)$
of an open subgroup of index $D$ in $\prod_{\ell \ne p} \SL_2(\Z_\ell)$
containing $\{\pm 1\}$,
then $\gamma_{\Fbar_p}(X_G) \ge \Phi_p(D)$.
\end{corollary}

\begin{proof}
Combine Proposition~\ref{P:special modular}
and Theorem~\ref{T:change}\eqref{I:simple square root}.
\end{proof}

\begin{remark}
\label{R:Atkin-Lehner}
One could also consider Atkin-Lehner quotients of $X_0(N)$
for $N$ prime to $p$.
These are generally not of the form $X_G$.
Nevertheless, 
gonality bounds tending to infinity for fixed $p$ can be obtained:
first apply Proposition~\ref{P:special modular} to get
a lower bound on $\gamma_{\F_{p^2}}(X_0(N))$,
next use Proposition~\ref{P:facts}(vi) to get a lower bound 
on the $\F_{p^2}$-gonality of any quotient of $X_0(N)$,
and finally apply Theorem~\ref{T:change}.
This works since the size of the Atkin-Lehner group (some power of $2$)
is asymptotically small compared to the index of 
the congruence subgroup $\Gamma_0(N)$ in $\PSL_2(\Z)$.
\end{remark}

\section{Level a power of $p$}
\label{S:power of p}

The necessary lower bound on the gonality of $\Ig(p^e)$ 
has been proved already by Schweizer, in a strong form:

\begin{theorem}[\cite{Schweizer2005}*{Lemma~1.5(d,e)}] \hfill
\label{T:Schweizer}
\begin{romanenum}
\item \label{I:1}
If $p\ge 7$, then
$\frac{p+13}{24} \le \gamma_{\Fbar_p}(\Ig(p)) \le \frac{p-1}{6}$.
\item \label{I:2}
If $e>1$ and $p^e \notin \{25,9,8,4\}$, then 
$\gamma_{\Fbar_p}(\Ig(p^e)) = p \gamma_{\Fbar_p}(\Ig(p^{e-1}))$.
\end{romanenum}
\end{theorem}

In fact, \cite{Schweizer2005} proves many more results.
The above are more than we need to deduce the following.

\begin{corollary}
\label{C:Psi}
Let $G$ be the inverse image under $S \surjects \Z_p^\times$
of an open subgroup of index $D$ in $\Z_p^\times$ containing $\{\pm 1\}$.
Then $\gamma_{\Fbar_p}(X_G) > D/12$.
\end{corollary}

\begin{proof}
First suppose that $X_G=\Ig(p^e)$ for some prime power $p^e>2$.
Then $D=p^{e-1}(p-1)/2$.
For $p^e \le 25$, we have $D<12$,
and $\gamma_{\Fbar_p}(X_G) \ge 1 > D/12$ trivially.
For $p>25$, Theorem~\ref{T:Schweizer}\eqref{I:1} gives
\[
\gamma_{\Fbar_p}(\Ig(p)) \ge \frac{p+13}{24} > \frac{p-1}{24} = \frac{D}{12}.
\]
For $p^e>25$ with $e>1$, we use induction on $e$,
with Theorem~\ref{T:Schweizer}\eqref{I:2} giving the inductive step.

Any other $X_G$ in Corollary~\ref{C:Psi} is a quotient of $\Ig(p^e)$
for some $p^e>2$, and the inequality for $X_G$ follows from
the inequality for $\Ig(p^e)$,
by Proposition~\ref{P:facts}(vi).
\end{proof}

\begin{remark}
An alternative approach to lower bounds on the gonality of Igusa curves
is to show that they have many points over certain finite fields,
to deduce that the gonality over these finite fields is large,
and then to apply Theorem~\ref{T:change}.
One can no longer use supersingular points, however, since these
are totally ramified in $\Ig(p^e) \stackrel{j}\to \PP^1$,
and hence their number does not grow with $e$.
Instead we could use {\em ordinary} points:
it follows from \cite{Pacheco1996}*{Corollary~2.13}
and Hurwitz class number estimates that $\#\Ig(q)(\F_q) \ge q^{1/2+o(1)}$
as $q \to \infty$.
Or one could use cusps, 
as in the proof of \cite{Schweizer2004}*{Theorem~6.1}, 
to get lower bounds on $\#\Ig(q)(\F_p)$, 
since the cusps split completely 
in $\Ig(p^e) \stackrel{j}\to \PP^1$ \cite{Katz-Mazur1985}*{Corollary~12.7.2}.
But the lower bounds on gonality obtained by these methods
are weaker than the ones we took from \cite{Schweizer2005}.
\end{remark}

\section{Group theory}
\label{S:group theory}

Here we study the open subgroups of $S$.
Let $S_p=\Z_p^\times$ and $S_{\ne p} = \prod_{\ell \ne p} \SL_2(\Z_\ell)$,
so $S = S_p \times S_{\ne p}$.
Given an open subgroup $G \le S$,
let $G_p$ and $G_{\ne p}$ be the images of $G$ in $S_p$ and $S_{\ne p}$,
respectively.

\begin{lemma}
\label{L:contains SL_2}
Let $B \in \Z_{>0}$.
Any open subgroup $H$ of $S_{\ne p}$ of index $\le B$ contains 
\[
        \prod_{\ell \le B!, \; \ell \ne p} \{1\} \times \prod_{\ell > B!, \; \ell \ne p} \SL_2(\Z_\ell).
\]
\end{lemma}

\begin{proof}
For each $\ell \ne p$, identify $\SL_2(\Z_\ell)$ with a
subgroup of $S_{\ne p}$ in the obvious way.
It suffices to show that $H$ contains $\SL_2(\Z_\ell)$
for each $\ell>B!$ with $\ell \ne p$.

The kernel of the action of $S_{\ne p}$ on the coset space $S_{\ne p}/H$
is a normal open subgroup $N \trianglelefteq S_{\ne p}$
contained in $H$.
Let $n:=(S_{\ne p}:N)$, so $n \le B!$.
Now $\ell>B! \ge n$, so $1/n \in \Z_\ell$, and 
\[
        \begin{pmatrix} 1 & 1 \\ 0 & 1 \end{pmatrix} 
        = \begin{pmatrix} 1 & 1/n \\ 0 & 1 \end{pmatrix}^n \in N,
\]
where the matrices belong to $\SL_2(\Z_\ell) \le S_{\ne p}$.
Similarly $\begin{pmatrix} 1 & 0 \\ 1 & 1 \end{pmatrix} \in N$.
But these two matrices generate the dense subgroup $\SL_2(\Z)$ 
of $\SL_2(\Z_\ell)$, so $\SL_2(\Z_\ell) \le N \le H$.
\end{proof}

\begin{lemma}
\label{L:open subgroups of S}
For each $B>0$,
there are at most finitely many open subgroups $G$ of $S$
such that $(S_p:G_p)<B$ and $(S_{\ne p}:G_{\ne p})<B$.
\end{lemma}

\begin{proof}
Fix $B$.
By Lemma~\ref{L:contains SL_2}, it suffices to consider instead
the situation in which $S_{\ne p}$ is replaced
by $S_L:=\prod_{\ell \in L} \SL_2(\Z_\ell)$
for a {\em finite} set $L$ of primes $\ne p$:
i.e., $G$ is now an open subgroup of $S_p \times S_L$,
$G_L$ is the image of $G$ in $S_L$,
and we are given $(S_p:G_p)<B$ and $(S_L:G_L)<B$.

Since $S_p$ and $S_{\ne p}$ are topologically finitely generated,
there are finitely many possibilities for $G_p$ and $G_L$.
Goursat's Lemma \cite{LangAlgebra}*{p.~75}
states that each possible $G$ is the inverse image under 
\[
        G_p \times G_L \surjects \frac{G_p}{H_p} \times \frac{G_L}{H_L}
\]
of the graph of an isomorphism 
\[
        \frac{G_p}{H_p} \to \frac{G_L}{H_L}
\]
for some normal open subgroups 
$H_p \trianglelefteq G_p$ and $H_L \trianglelefteq G_L$.
By the finite generation again,
it suffices to bound $(G_p:H_p)=(G_L:H_L)$.
It is bounded by the supernatural number $\gcd(\#S_p,\#S_L)$,
which is finite, since $S_p$ has a pro-$p$ open subgroup,
while $S_L$ has an open subgroup of order prime to $p$.
(See \cite{SerreGaloisCohomology}*{I.\S1.3} for the notion
of supernatural number.)
\end{proof}

\section{The general case of Theorem~\ref{T:main}}
\label{S:general case}

\begin{proof}[Proof of Theorem~\ref{T:main}]
Let $m>0$.
We will show that $\gamma_{\Fbar_q}(X_{G_i})>m$
for all but finitely many $i$.
Let $\Phi_p$ be as in Corollary~\ref{C:Phi}.
Choose $B_0$ such that $\min\{\Phi_p(B),B/12 \} > m$ for all $B \ge B_0$.
By Lemma~\ref{L:open subgroups of S},
all but finitely many $G_i$ in our sequence
have either $(S_p:(G_i)_p) \ge B$ or $(S_{\ne p}:(G_i)_{\ne p}) \ge B$.
For each such $i$, $X_{G_i}$ dominates an $X_G$ 
with $G$ as in Corollary \ref{C:Phi} or Corollary \ref{C:Psi};
then $\gamma_{\Fbar_p}(X_G)$ exceeds either $\Phi_p(B)$ or $B/12$.
By Proposition~\ref{P:facts}(vii),
\[
        \gamma_{\Fbar_p}(X_{G_i}) \ge \gamma_{\Fbar_p}(X_G) 
	\ge \min\{ \Phi_p(B), B/12 \} > m.
\]
\end{proof}

\section{Image of Galois}
\label{S:uniform boundedness}

\begin{proof}[Proof of Theorem~\ref{T:uniform boundedness}]
We may assume that $k$ is algebraically closed.
By Theorem~\ref{T:main},
there are only finitely many subgroups $G \le S$
containing $\{\pm1\}$ such that $\gamma_{\Fbar_p}(X_G) \le d$.
Choose $N_{p,d}$ such that $N_{p,d} \ge 2(S:G)$ for every such $G$.

Let $K$ be a field of degree $\le d$ over $k(t)$,
and let $E$ be an elliptic curve over $K$
with $j(E)$ not algebraic over $k$ (i.e., not in $k$).
Write $K=k(C)$, where $C$ is a curve over $k$ with $\gamma_k(C) \le d$.
Define $H := \rho_E(\Gal(K^s/K))$.
Since $k$ is algebraically closed, $H \subseteq S$.
We want $(S:H) \le N_{p,d}$.

Suppose not.
If $(S:H)$ is infinite,
then since $S/H$ is a profinite group,
we can find a group $H'$ with $H \le H' \le S$
and $N_{p,d} <  (S:H') < \infty$.
If $(S:H)$ is finite, let $H'=H$.
In either case, let $H''$ be the group generated by $H'$ and $-1$,
so $N_{p,d}/2 < (S:H'') < \infty$.
By definition of $N_{p,d}$, the group $H''$ does not equal any of the
groups $G$, so $\gamma_{\Fbar_p}(X_{H''})>d$.
Equivalently, by Proposition~\ref{P:facts}(ii), $\gamma_k(X_{H''})>d$.

The curve $X_{H''}$ is defined as a quotient of some $X(p^e;N)$.
Choosing level structure for $E$ over $K^s$ gives a point 
in $X(p^e;N)(K^s)$,
and the action of $\Gal(K^s/K)$ moves this point within the $H''$-orbit,
since $H \subseteq H''$,
so the image point in $X_{H''}(K^s)$ is $K$-rational.
This point in $X_{H''}(K)$ may be viewed as a rational map 
$C \dashrightarrow X_{H''}$,
and this map is non-constant since the composition
$C \dashrightarrow X_{H''} \stackrel{j}\to X(1) \isom \PP^1$
corresponds to $j(E) \in K - k$.
Proposition~\ref{P:facts}(vii) implies 
$\gamma_k(X_{H''}) \le \gamma_k(C) \le d$,
contradicting the previous paragraph.
\end{proof}

\section*{Acknowledgements} 

Most of all I thank Andreas Schweizer for several discussions,
for suggesting several references, and for pointing out that existing
versions of the Castelnuovo-Severi inequality seem to require a
perfect field of constants.
I also thank Matt Baker, from whom I first learned
the idea in \cite{Harris-Silverman1991}
that the Castelnuovo-Severi inequality 
could be used to bound change in gonality under algebraic extensions.
I thank Brian Osserman for pointing out that the lower bound
on the geometric gonality of the general curve of genus $g$ 
in characteristic $p$ could be deduced from \cite{Osserman2005}.
I thank Doug Ulmer for a remark about cusps of Igusa curves.
Finally I thank Alina Cojocaru, Chris Hall, Dinesh Thakur, and Doug Ulmer
for asking questions that inspired this paper.

\end{document}